\documentclass[microtype,psamsfont]{amsart}

\usepackage{graphicx}
\usepackage[dvips]{epsfig}
\usepackage{psfrag}
\usepackage{pinlabel}
\usepackage{amsmath}
\usepackage{amsfonts}
\usepackage{latexsym}
\usepackage{amssymb}
\usepackage[usenames]{color}
\usepackage{amsthm}
\input xy
\xyoption{all}
\author{Baptiste Chantraine }
\address{Universit\'e du Qu\'ebec \`a Montr\'eal, Montr\'eal, Canada, H3C 3P8}
\address{Universit\'e Libre de Bruxelles, Bruxelles, Belgique, 1050}
\email{bchantra@ulb.ac.be}

\newcommand{\difff}{\omega}
\newcommand{\del}{\partial}
\newcommand{\al}{\alpha}
\newcommand{\eps}{\epsilon}
\newcommand{\emb}{\hookrightarrow}

\newcommand{\appl}{\rightarrow}

\theoremstyle{remark}
\newtheorem{defn}{Definition}[section]
\newtheorem{Rem}[defn]{Remark}
\newtheorem{Conj}[defn]{Conjecture}
\theoremstyle{plain}
\newtheorem{Thm}{Theorem}[section]
\newtheorem{Prop}{Proposition}[section]

\newtheorem{Cor}[defn]{Corollary}
\makeatletter
  \let\c@Thm=\c@defn
\makeatother
\makeatletter
  \let\c@Prop=\c@defn
\makeatother

\def\Ex{{\bf Example. }}
\def\co{\colon\thinspace}

\DeclareMathAlphabet{\mathdj}{U}{msb}{m}{n}
\newcommand{\R}{\ensuremath{\mathdj{R}}}
\newcommand{\Z}{\ensuremath{\mathdj{Z}}}
\newcommand{\Qu}{\ensuremath{\mathdj{H}}}
\newcommand{\C}{\ensuremath{\mathdj{C}}}

\def\L{\mathcal{L}}

\begin{document}
\title{Lagrangian concordance of Legendrian knots }

\begin{abstract}In this article, we define the notion of a Lagrangian concordance between two Legendrian knots analogous to smooth concordance in the Legendrian context. We show that Legendrian isotopic Legendrian knots are Lagrangian concordant. The focus is primarily  on the algebraic aspects of the problem. We study the behavior of the classical invariants (namely the Thurston--Bennequin number and the rotation number) under this relation, and provide some examples of nontrivial Legendrian knots bounding Lagrangian surfaces in $D^4$.
\end{abstract}
\thispagestyle{empty}
\maketitle


\section{Introduction}
A contact structure $\xi$ on a $3$--manifold $Y$ is a completely non-integrable $2$--plane field (ie it is locally defined as
the kernel of a $1$--form $\al$ such that $\al\wedge d\al\not= 0$). Through this paper $Y$ will be oriented and the
contact structure will be assumed to be transversally orientable and positive. Transversally orientable  means that the $1$--form $\al$ can be globally defined and the positivity condition means that $\al\wedge d\al >0$ (note that this is really a condition on the contact structure and not on the $1$--form defining it).

The completely non-integrable condition implies that
any surface embedded in $Y$ cannot be tangent to $\xi$. However many one-dimensional submanifolds tangent to $\xi$ exist. Such submanifolds are called {\it Legendrian}. A Legendrian knot in $(Y,\xi)$ is a map $\gamma \co S^1\emb Y$ such that for all
$s\in S^1,T_s\gamma(S^1)\subset\xi$ (this is equivalent to $\gamma^*\al=0$). We usually denote by $K$ the image of $\gamma$ and call the image Legendrian knot also. Two Legendrian knots are \textit{Legendrian isotopic} if there is a (smooth)
isotopy $I$ between them such that  $I(\cdot ,t)$ is a Legendrian knot for all $t$.

The symplectisation of a contact manifold $(Y,\xi)$ is the manifold $(\R\times Y, d(e^v\al))$. Most of contact geometry of $Y$ can reformulated in terms of equivariant (or invariant) symplectic geometry of $\R\times Y$.

The goal of this paper is to define relations on the set of Legendrian knots, namely \textit{Lagrangian cobordism} and \textit{Lagrangian concordance}, and to study their properties.
\begin{defn}\label{lagcob}
 Let $\Sigma$ be a compact orientable surface with two points $p^-$ and $p^+$ removed. Let $(u^-,s) \in (-\infty,-T)\times S^1$ and $(u^+,s)\in (T,+\infty)\times S^1$ be cylindrical coordinates around $p^-$ and $p^+$ respectively. Let $\gamma^-$ and $\gamma^+$ be Legendrian knots in $Y$. Then:
 \begin{enumerate}
 \item {$\gamma^-$ is  {\it Lagrangian cobordant}  to $\gamma^+$ (written $\gamma^-\prec_\Sigma\gamma^+$) if there exists a Lagrangian embedding
$$L\co \Sigma\emb \R\times Y$$
such that:
\begin{itemize}
 \item{$L(u^-,s)=(u^-,\gamma^-(s))$}
\item{$L(u^+,s)=(u^+,\gamma^+(s))$}
\end{itemize}}
\item {$\gamma^-$ is Lagrangian concordant to $\gamma^+$ if $\gamma^-\prec_\Sigma\gamma^+$ and $\Sigma$ is a cylinder (we will denote this particular case by $\gamma^-\prec\gamma^+$)}.
 \end{enumerate}
\end{defn}
In this paper we will mostly be concerned with Lagrangian concordance.

Our first result is the following Theorem.
\begin{Thm}\label{Theoremprinc}

Let $K^-$ and $K^+$ be two Legendrian knots in $Y$ and let $I\co S^1\times [0,1] \appl Y$ be a Legendrian isotopy between $K^-$ and $K^+$. Then there exists a Lagrangian cylinder $C$ such that $K^-\prec_C K^+$.
\end{Thm}
In other words, the Lagrangian concordance relation descends to the Legendrian isotopy classes of Legendrian knots. \ref{Theoremprinc} will be proved in \ref{isotopieconc}.

The next theorem gives the relationship between the classical invariants of two Legendrian knots which are Lagrangian cobordant.
\begin{Thm}\label{cobinvariance}
 If $\gamma^-\prec_\Sigma \gamma^+$ then:
\begin{align*}
 &r(\gamma^-,[S])=r(\gamma^+,[S\cup \Sigma])\\
&tb(\gamma^+)-tb(\gamma^-)=-\chi(\Sigma)
\end{align*}
Where $tb$ refers to the Thurston--Bennequin number, $r$ the rotation number.
\end{Thm}
We give a proof of \ref{cobinvariance} in \ref{Behaviour}.

Combining \ref{cobinvariance} and a Bennequin type inequality we prove, in \ref{cob}, the following theorem:
\begin{Thm}\label{lagsurface}
Let $\Sigma$ be an oriented Lagrangian submanifold of a Stein surface such that $\del\Sigma$ is a Legendrian submanifold $K$
of $\del X$. Then the following holds:
$$r(K,[\Sigma])=0$$
$$tb(K)=2g(\Sigma)-1=TB(\L(K))$$
Furthermore, if $X\simeq D^4$ with its standard Stein structure, then $$g(\Sigma)=g_s(\L(K)),$$
where $\L(K)$ is the smooth isotopy type of $K$ , $TB(\L(K))$ is the maximal Thurston--Bennequin number of Legendrian representatives of $\L(K)$ and $g_s(\L(K))$ is the $4$--ball genus
of $\L(K)$.
\end{Thm}

This result is the principal topological motivation for our definitions. It provides a criterion to compute the $4$--ball genus of some knots. In \ref{algknot} we provide a large class of examples where \ref{lagsurface} applies.

We then conclude in \ref{conclusion} with some remarks and potential applications.
\subsection*{Addendum:}
Since the appearance of the first version of this paper, new developments are in preparation. Using a forthcoming result of Ekholm, Honda and K\'alm\'an \cite{KalHoEk} we can show that there exists concordances which do not come from Legendrian isotopy. For instance, using their proof of the existence of a Lagrangian cobordism realizing the $1$--smoothing of crossing in the Lagrangian projection, one can easily see that the $9_{46}$ knot in the Rolfsen table admits a Legendrian representative, with maximal Thurston--Bennequin invariant, which is the boundary of a Lagrangian disk. Hence the trivial Legendrian knot $K_0$ is Lagrangian concordant to this one.

In a forthcoming note we will show that this concordance is not reversible, meaning that this knot cannot be concordant to $K_0$. Thus Lagrangian concordance is indeed not a symmetric relation.
\subsection*{Acknowledgments:} This work makes up part of my Ph.D. thesis in Universit\'e du Qu\'ebec \`a Montr\'eal under the supervision of Olivier Collin, whom I warmly thank for many helpful comments and discussions. I also wish to thank Paolo Ghiggini and Tam\'as K\'alm\'an for inspiring conversations, and Tobias Ekholm for pointing out another way to see the behavior of the Thurston--Bennequin number under cobordism.

 \section{Preliminaries}
In this section we collect once and for all the basic facts of contact and symplectic geometry we apply through the paper. We first consider some generalities about the symplectisation of a contact manifold and its Hamiltonian diffeomorphism group (\ref{symplectisation}). Then in \ref{Clainv} we give the facts about Legendrian knots which will be needed later. We also give several ways to interpret the Thurston--Bennequin number and the rotation number in our context.

But we first fix the notation we will use through the paper.
\begin{itemize}\item {The letter $s$ will always refer to a parameter in $S^1=\R/2\pi\Z$.}
\item  {We assume that our cylinders are parametrized by $\{(u,s)| u\in\R
  ,\,s\in\R/2\pi\Z\}$.}
\item {The letter $v$ denotes the parametrization of $\R$ in the
  symplectisation $\R\times Y$.}
\item {The Reeb vector field of a contact form is the unique vector field
  $R_\al$satisfying: $\al(R_\al)=1$ and $d\al(R_\al,\cdot)=0$.}
\item {Since $\xi$ is a symplectic vector bundle (with the form $d\al$), it
  admits a compatible almost-complex structure (we always denote by
  $J$). It extends to $T(\R\times Y)$ by setting $$J\frac{\del}{\del
    v}=R_\al.$$}
\item {A symplectic bundle with a compatible almost complex structure is
  naturally equipped with a Euclidian (resp. Hermitian) metric given
  by $g(V,W)=\difff(JV,W)$. Under this context, we assume that all
  trivializations are orthonormal (resp Hermitian).}
\item {We let $\langle\cdots\rangle_\C={\text span}_\C(\cdots)$ be the
  complex span of the given vectors and $\langle\cdots\rangle={\text
    span}_\R(\cdots)$ the real span.}
\item {We consider that $(1,i,j,k)$ is the real basis of the quaternionic
  line $\Qu$.}
\item {The standard contact structure $\xi_0$ on $S^3=\del
  D^4\subset\C^2(\simeq\Qu)$ will be the one defined by the complex
  tangencies $\xi_{0,p}=T_pS^3\cap i(T_pS^3)$. However, in order to
  simplify the notation, we will sometimes use $\eta_{0,p}= T_pS^3\cap
  j(T_pS^3)$. We will refer the first description as the $i$--convex
  contact structure and the second as the $j$--convex one. These two
  contact structures are obviously contactomorphic by a linear
  transformation in $\R^4$.}
\item {In the $i$--convex sphere, we usually denote by $K_0$ the trivial
  Legendrian knot given by $\gamma_0(s)=(\cos s,0,\sin s,0)$. Its
  analogue in the $j$--convex sphere is $\gamma_0(s)=(\cos s,\sin
  s,0,0)$.}
\end{itemize}
\subsection{Symplectisation}\label{symplectisation}
In this section we study some particular examples of Hamiltonian diffeomorphisms of the symplectisation of a contact manifolds. \ref{Cont-Ham} will be the main tool to prove  \ref{Theoremprinc}.

Let $(Y,\xi)$ be a manifold together with an hyperplane distribution $\xi$.
Let $X=\{(p,\lambda)\vert \ker\lambda=\xi_p\}\subset T^*Y$ be the annihilator of this distribution minus the $0$--section.
Consider the Liouville form $\theta$ of $T^*Y$ defined as $\theta_{(p,\lambda)}(V)=\lambda_p(\pi_*(V))$. If we also denote its restriction to $X$ by $\theta$ then we have the following obvious proposition.
\begin{Prop}
$\xi$ is a contact structure if and only if $(X,d\theta)$ is a symplectic manifold.
\end{Prop}

Under these hypothesis $(X,d\theta)$ is called the symplectisation of $(Y,\xi)$. This definition coincides with the previous one. Indeed a contact form is nothing but a section of $X$ as an $\R^*$--bundle over $Y$ and hence gives an isomorphism of $(X, d\theta)$ with $(\R^*\times Y, d(w\al))$. The positive connected component is then identified with $(\R\times Y,d(e^v\al))$.

We denote by $\text{Symp}(X,\difff)$ the set of symplectomorphisms of $(X,\difff)$ and $\text{Ham}(X,\difff)$ the set of Hamiltonian diffeomorphisms. Also, we denote by $\text{Cont}(Y,\xi)$ the set contactomorphisms of $(Y,\xi)$. We denote the formal Lie algebra of those groups by $\mathfrak{S}(X,\difff)$ (the time-dependent {\it symplectic vector fields}), $\mathfrak{H}(X,\difff)$ (the time-dependent {\it Hamiltonian vector fields}) and $\mathfrak{C}(Y,\xi)$ (the time-dependent {\it contact vector fields}) respectively. More precisely:
\begin{align*}
 \mathfrak{S}(X,\difff)&=\{W_t\in\Gamma(TX)\vert \L_{W_t}\difff=0\}\\
 \mathfrak{H}(X,\difff)&=\{W_t\in\Gamma(TX)\vert \difff(W_t,\cdot)=dH_t\}\\
\text{and }\mathfrak{C}(Y,\xi)&=\{V_t\in\Gamma(TY)\vert \L_{V_t}\al=R_\al(\al(V_t))\al\}\\
\end{align*}

The choice of a contact form for a contact structure $\xi$ on $Y$ yields an isomorphism between the algebra of time-dependent functions $H\co Y\times \R/\Z\appl\R$ and contact vector fields by the following correspondence:
\begin{align}\label{hamiltonniendecontact}
\al(V_t)=H(\cdot,t),\\
dH_t=dH_t(R_\al)\al-V_t\iota d\al,
\end{align}
where $\iota$ denotes the inner product $\iota\co \Gamma(TY)\otimes\Omega^2(Y)\mapsto \Omega^1(Y)$.
Those two equations uniquely determine a time-dependent contact vector field given a function $H$ and reciprocally a time-dependent contact vector field defines a {\it contact Hamiltonian} function by the first equation.

Recall that any diffeomorphism $f$ of $Y$ induces a symplectomorphism $F$ of $T^*Y$ by the formula $F(p,\lambda)=(f(p),(f^{-1})^*(\lambda))$. Moreover $f$ is a contactomorphism if and only $F(X)=X$.
This implies that we have a natural homomorphism
\[ \begin{array}{cccc}
i\co & \text{Cont}(Y,\xi) & \appl & \text{Symp}(X,\difff) \\
     & f                &\appl & F\vert_X
\end{array} \]
which induces a homomorphism $$i_*\co \mathfrak{C}(Y,\xi)\appl \mathfrak{S}(X,\difff)$$ (if $Y$ is non-compact this homomorphism is still defined using local flow). With those notations we have the following:

\begin{Prop}\label{Cont-Ham}
$\text{im}(i_*)\subset \mathfrak{H}(X,\difff)$ and thus $$i\co \text{Cont}(Y,\xi) \appl \text{Ham}(X,\difff)$$
ie the lift of a contactomorphism of $Y$ to a symplectomorphism of the symplectisation $X$ is a Hamiltonian diffeomorphism.
\end{Prop}
\begin{proof}
The symplectomorphism defined above is actually exact:
\begin{align*}
 F^*(\theta)_{(p,\lambda)}(W)& =\theta_{(f(p),(f^{-1})^*\lambda)}(F_*(W))\\
& =(f^{-1})^*\lambda_{f(p)}(\pi_*\circ F_*(W))\\
& =\lambda_p(f^{-1}_*\circ f_*\circ\pi_* (W))=\lambda_p(\pi_*(W))\\
& =\theta_{(p,\lambda)}(W)
\end{align*}
For a one parameter family of contactomorphisms $f_t$ generated by a time-dependent vector field $V_t$ this gives:
\begin{align}\label{hamiltonien}
F_t^*(\theta)=\theta
\end{align}
If one sets $\widetilde{V_{t_0}}_p=\frac{d}{dt}\vert_{t=t_0}F_t(p)$ and differentiates \ref{hamiltonien} with respect to $t$ we get

\begin{align*}
0=\L_{\widetilde{V_t}}\theta&=d(\theta(\widetilde{V_t}))+\widetilde{V_t}\iota d\theta\\
\Rightarrow \widetilde{V_t}\iota \difff&=-d(\theta(\widetilde{V_t}))
\end{align*}
which proves that $\widetilde{V_t}$ is Hamiltonian.\\
By construction $\widetilde{V_t}=i_*(V_t)$.
\end{proof}

Hence we have the following commutative diagram
 \begin{align*}
\xymatrix{ {\mathfrak{C}(Y,\xi)} \ar[d]_{exp} \ar[r]^{i_*} & {\mathfrak{H}(X,\difff)} \ar[d]_{exp} \\
                    {\text{Cont}(Y,\xi)} \ar[r]^{i}                  & {\text{Ham}(X,\difff)}                                          }
                    \end{align*}
(where we restrict $exp$ to its domain of definition).

If one writes the contact vector field $V$ as a sum $HR_{\al}+W$ where $W\in \xi$ (thus $H$ is the contact Hamiltonian) then the Hamiltonian function generating the lifts of the contactomorphism is given by $\theta(V)$ and thus by the function $e^t\cdot H$. A simple computation then gives that the Hamiltonian vector field is given by: $$-dH_t(R_{\al})\frac{\del}{\del t}+V_t$$

\subsection{Legendrian knots}\label{Clainv}

We discuss here the facts that will be needed later about Legendrian knots.

By a Legendrian isotopy from $\gamma_0$ to $\gamma_1$ we mean a smooth isotopy $$I\co S^1\times [0,1]\appl (Y,\xi)$$ such that:\\
\begin{itemize}
\item{for all $t_0\in [0,1]$, $I(\cdot,t_0)\co S^1\appl (Y,\xi)$  is a Legendrian knot.}
\item{$I(s,i)=\gamma_i(s),i=0,1$.}
\end{itemize}

It is well known (see Geiges \cite[Theorem 2.41]{gei} for a proof) that a Legendrian isotopy can be realized by an ambient contact isotopy, ie
there exists a smooth family of compactly supported contactomorphisms $f_t$ such that $I(s,t_0)=f_{t_0}(\gamma_0(s))$.
So any Legendrian isotopy is realized as the evolution of the original knot along a contact flow, the ``Lie algebra'' of Legendrian isotopies is then the contact vector fields along a Legendrian knot (note here that this is not the case for Lagrangian isotopy where we need to consider conformally symplectic fields).

To a null-homologous Legendrian knot one can associate three classical invariants:
\begin{itemize}
 \item {The smooth type of the knot}
\item{its Thurston--Bennequin invariant}
\item{its rotation number}
\end{itemize}

Since in the proofs of the results in the present paper we used several different descriptions of those invariants, we describe them here.

\paragraph{The Thurston--Bennequin invariant.}

To define this invariant we need an orientation of the knot, the resulting invariant will not depend on this choice. We assume that $K$ is oriented by its parametrization $\gamma$. We denote by $\tau$ its unit positive tangent vector.

The Thurston--Bennequin invariant $tb_K$ is the homotopy class of the trivialization of the tubular neighborhood of $K$ given as follows: first take a vector field $V$ along $K$ inside $\xi$ transverse to $\tau$ so that $(\tau,V)$ gives the positive orientation of $\xi$ and take a vector field along $K$ transverse to contact structure giving the positive transverse orientation (for instance we can take $tb_K=(J(\tau),R_\al)$).

If $K$ is null-homologous we can choose a Seifert surface $S$ for $K$. The surface $S$ gives an orthogonal trivialization of $\mathcal{N}(K)$ and the difference between this one and $tb$ is thus an element of $\pi_1(SO(2))=\Z$. We can therefore assign to $tb$ an integer (independent of $S$) which we denote by $tb(K)$. The number $tb(K)$ is therefore computed as the intersection number between a small push off of $K$ along $R_\al$ and $S$.

Note also that $(K,tb_K)$ is a framed submanifold of $Y$ of codimension 2. Via the Thom--Pontryagin construction, it corresponds to a map $h\co Y\appl S^2$. Recall that this map is defined as follows: a framing identifies a neighborhood of $K$ with $S^1\times D^2$. The map from $Y$ to $S^2\simeq D^2\cup\{*\}$ projects $S^1\times D^2$ to $D^2$ and sends everything else to $*$. The choice of $S$ gives a (framed) cobordism $S_0$ from $K$ to the trivial knot (with a framing $\sigma$ obtained by extending the trivialization $tb_K$ to $S_0$), so there is a map $h'$ homotopic to $h$ which factors through the fundamental map $[Y]\co Y\appl S^3$.
\begin{align*}
 \xymatrix{ {Y} \ar[r]^{[Y]} \ar[d]_{h'} & {S^3} \ar[dl] \\
            {S^2}
}
\end{align*}

Consequently the homotopy type of $h'$ is characterized by a homotopy class of maps from $S^3\appl S^2$. Recall that the free homotopy group $[S^3,S^2]$ is isomorphic to $\Z$ (generated by the Hopf fibration), therefore $h'$ is characterized by this integer and we claim the following:

 \begin{Prop}\label{thompontTB} Let $n\in \Z$ be the integer associated to $tb_K$ and $[S]$ as above then $tb(K)=-n$.
 \end{Prop}

\begin{proof}
 We first notice that in the case of the trivial Legendrian knot $K_0=\{z_2=0\}$ in $S^3=\{(z_1,z_2)\in\C^2\simeq\Qu| |z_1|^2+|z_2|^2=1\}$ with the $j$--convex contact structure, then $tb_{K_0}$ is $( (ie^{i\theta},0),(je^{i\theta},0))$ which is the trivialization induced by the Hopf map. In this particular case $tb(K_0)=-1$ (by a direct computation or by \ref{cob}). So the proposition is verified in that case.
 
Consider now the general case where $tb(K)=n$. Choose a cobordism from $K$ to $K_0$ ($K_0$ in a standard Darboux chart). Since $tb(K)=n$, the extension of the trivialization to the cobordism gives a trivialization on $K_0$ which is $-n$ times the one given by $tb_{K_0}$. The Thom--Pontryagin map is thus $-n$ times the Thom--Pontryagin map associated to $tb(K_0)$, which completes the proof. 
\end{proof}
 \begin{Cor}
 The Thurston--Bennequin number of a null-homologous Legendrian knot in $Y$ is completely determined by the Thom map associated to the trivialization given by $tb(K)$.
 \end{Cor}
\begin{Rem}
 If $K$ happens to be non null-homologous we still can associate to a framing an integer which is well defined modulo twice the divisibility of $[K]\in H_1(Y)/\text{tor}H^1(Y)$.  \ref{thompontTB} holds also in that case.
\end{Rem}

\paragraph{The rotation number.}

This number depends both on the orientation of $K$ and on the homology class of a Seifert surface $S$. The rotation number is the relative Euler class of $\xi\vert_S$ with respect to $\tau$, $r(K,[S])=e(\xi\vert_S,\tau)\in H^2(S,\del S)$. The Euler class can, in fact, be replaced by the first  Chern class since $d\al\vert_\xi$
is a symplectic form on $\xi$.

The relative Chern class can be computed as follows:
consider the trivial complex line bundle $D^2\times\C\appl D^2$ and use the trivialization of the complex line bundle $\xi|_{K}$ given by $\tau(s)=\frac{\del}{\del s}$ to construct a complex line bundle $\xi'$ over $S':=S\cup D^2$ (identifying $\tau$ along $K$ with $1\in\C$ along $S^1=\del D^2$). We set $c_1(\xi\vert_S,\tau)=c_1(\xi')$ using the isomorphism $H^2(S,\del S)\simeq H^2(S')\simeq\Z$.

 The rotation number can be computed noticing that $\xi\vert_S$ is a complex line bundle over $S$ which has the homotopy type of a wedge of circles, hence $\xi\vert_S$ is trivializable over $S$. Denote such an Hermitian trivialization by $\sigma\co \xi\vert_S\simeq S\times\C$  (let say given by sections $V$ and $JV$). Then $\sigma\vert_{K}$ is a trivialization of $\xi$ along $K$ (notice that since $K$ represents a product of commutators in $\pi_1(S,x_0)$ and since $\pi_1(U(1),Id)$ is Abelian, the homotopy type of this trivialization does not depend on $\sigma$). With respect to this trivialization, $\tau$ becomes a loop in $U(1)$ and the homotopy type of this loop in $\pi_1(U(1),Id)\simeq\Z$ is exactly the rotation number.

\Ex
Consider $S^3=\del D^4$ with its $i$--convex structure and $K\co [0,2\pi]\appl S^3$ a Legendrian knot. We have $$\xi_{(z_1,z_2)}=\langle -z_2\frac{\del}{\del z_1}+z_1\frac{\del}{\del z_2}\rangle_\C$$ hence $\xi_0$ is already trivialized over $S^3$. So we can choose $\sigma$ to be this trivialization restricted to $S$ . It follows that $$\tau(s)=e^{i\theta(s)}(-z_2\frac{\del}{\del z_1}+z_1\frac{\del}{\del z_2})$$ and therefore $$r(K_0)=\frac{1}{2\pi} \int\limits_{S^1}d\theta.$$

There is a third way to compute the rotation number in terms of symplectic geometry. Consider in the symplectisation of $Y$ the trivial cylinder $C=\R\times K$. At any point of $X=\R\times Y$ the (symplectic) tangent space splits as $$T_p X=\langle\xi_p\rangle\oplus\langle\frac{\del}{\del v},R_\al\rangle.$$ This splitting is complex as well as symplectic as it can easily be checked. So we may  also trivialize $TX\vert_S$ using $\sigma$. Consider $K$ as a loop on the Lagrangian cylinder $C$ so that one can interprets, via our trivialization, the loop $T_{K(s)}C$ as a loop of Lagrangian planes $l(s)$ in $\C^2$. One can consequently associate to this loop its Maslov index $\mu(K,S)=[l(s)]\in\pi_1(\Lambda(2),\Lambda_0)\simeq H^1(\Lambda(2))\simeq\Z$, where we denote by $\Lambda(2)$ the set of Lagrangian planes in $\C^2$. Again, note that this number does not depend on $\sigma$ because $K$ represents a product of commutator in $\pi_1(S,x_0)$.

A good discussion on the Maslov index can be found in Viterbo \cite[Section 1]{Vitind}. We discuss here the definition of the Maslov index in our particular case.

Let $x=(1,0)$ and $y=(0,1)$ in $\C^2$ be the image of $v$ and $\frac{\del}{\del t}$ by the trivialization $\sigma$,
since the loop of Lagrangian planes is given by $$T_{K(s)}C=\langle e^{i\theta(s)}v(s),\frac{\del}{\del t}\rangle$$ we have $l(s)=\langle e^{i\theta(s)}\cdot x,y\rangle$.

Denote by $\Lambda_0$ the Lagrangian plane in $\C^2$ given by $\langle x,y+iy\rangle$ (the image under the trivialization of $\langle V,\frac{\del}{\del t}+R_\al\rangle$). The Maslov index of the loop is then, by definition (see \cite{Vitind}), the intersection of $l(s)$ with the Maslov cycle $$\{\Lambda\in\C^2-\text{Lagrangian}\vert \Lambda\cap \Lambda_0\not=\{0\}\}.$$
These intersections arise when $\theta(s)= 0\,\text{mod}\,\pi$ and then the intersection form is given by: $$Q(s)(x_1,x_2)=\tan\theta(s)\cdot x_1^2+x_2^2.$$ Hence intersections are positive when $\theta '(s)>0$ and negative when $\theta '(s)<0$ (note that the transversality condition is achieved assuming $\theta '(s)\not=0$ at intersection points). So we obtain the formula $$\mu(l)=2\int\limits_{S^1}d\theta$$ which implies $\mu(l)=2r(K,S)$, so that $$\mu(K,S)=2r(K,[S]).$$


 \section{Lagrangian concordance and its relation to Legendrian isotopy}$\label{Princ}$
\subsection{Definition of Lagrangian concordance}
        Let $Y$ be a smooth oriented $3$--manifold, $\xi$ a positive contact structure and
        $(X,\difff)\cong (\R\times Y,d(e^v\al))$ be the symplectisation of $Y$.
We recall here the definition given in introduction.

\medskip
{\bf Definition \ref{lagcob}}
\textsl{Let $\Sigma$ be a compact orientable surface with two points $p^-$ and $p^+$ removed. Let $(u^-,s) \in (-\infty,-T)\times S^1$ and $(u^+,s)\in (T,+\infty)\times S^1$ be cylindrical coordinates around $p^-$ and $p^+$ respectively. Let $\gamma^-$ and $\gamma^+$ be Legendrian knots in $Y$. Then:
 \begin{enumerate}
 \item {$\gamma^-$ is  {\it Lagrangian cobordant}  to $\gamma^+$ (written $\gamma^-\prec_\Sigma\gamma^+$) if there exists a Lagrangian embedding
$$L\co \Sigma\emb \R\times Y$$
such that:
\begin{itemize}
 \item{$L(u^-,s)=(u^-,\gamma^-(s))$}
\item{$L(u^+,s)=(u^+,\gamma^+(s))$}
\end{itemize}}
\item {$\gamma^-$ is Lagrangian concordant to $\gamma^+$ if $\gamma^-\prec_\Sigma\gamma^+$ and $\Sigma$ is a cylinder (we will denote this particular case by $\gamma^-\prec\gamma^+$)}.
 \end{enumerate}}

\begin{figure}[ht!]
\labellist
\small\hair 2pt
\pinlabel {$\R\times M$} [tl] at 40 195
\pinlabel {$\gamma^-$} [br] at 19 148
\pinlabel {$C$} [br] at 132 189
\pinlabel {$\gamma^+$} [bl] at 318 185
\pinlabel {$-T$} [tl] at 62 2
\pinlabel {$T$} [tl] at 202 2
\pinlabel {$\nu$} [bl] at 350 5
\endlabellist
\begin{center}
\includegraphics[height=5cm]{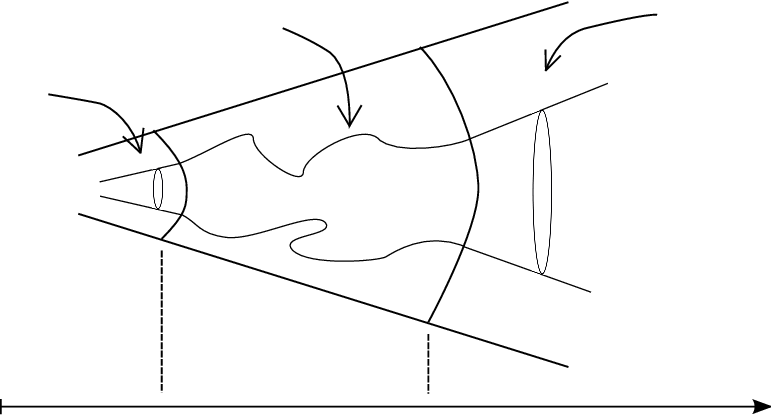}
 \end{center}
 \caption{A Lagrangian concordance.}
 \label{Concor}
 \end{figure}
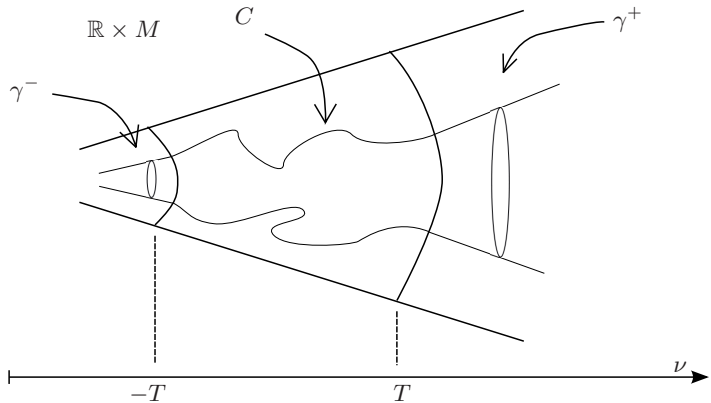

The previous definition is motivated by the following: any $\R$--invariant Lagrangian submanifold of $Y$ projects to a Legendrian submanifold of $Y$, and any Legendrian submanifold lifts to an $\R$--invariant Lagrangian submanifold. In particular any Legendrian submanifold is Lagrangian concordant to itself.

However, for the theory of Lagrangian concordance to be more intimately related to Legendrian knot theory, we want a relation up to Legendrian isotopy rather than a relation on the Legendrian submanifolds themselves. For this we have to prove that any Legendrian isotopy in the contact manifold $Y$ gives rise to a Lagrangian cylinder in the symplectisation $X$. The proof of this fact is the main goal of this section. But first some remarks about the differences between the Lagrangian concordance and the topological (smooth) concordance are in order.
\begin{Rem}
In the case of a smooth concordance, the condition of symmetry needed to obtain an equivalence relation is automatic. One simply needs to reverse the $u$ parameter in the embedding defining the concordance. In our case, the diffeomorphism $\phi\co  X\appl X$ sending $(x,v)$ to $(x,-v)$ is not a symplectomorphism, since $\phi^*(d(e^v\al))=d(e^{-v}\al)$ and $\phi$ sends parts of big volume to parts of small volume. Therefore it is not automatic that $K^-$ is Lagrangian concordant to $K^+$ implies that $K^+$ is Lagrangian concordant to $K^-$. This is actually not true as mentioned in the Addendum in Introduction.
\end{Rem}
\begin{Rem}

         On the symplectisation $X$, the vector field $\frac{\del}{\del v}$ expands the symplectic form $\omega$, so that any translation of a Lagrangian submanifold along this vector field will remain Lagrangian. This allows us to prove that $\prec$ is transitive. This relation thus resembles a partial order on set of the isotopy classes of Legendrian submanifolds. Whether or not it is a real partial order is not clear yet and seems to be out of reach since one must prove that $\gamma^-\prec\gamma^+$ and $\gamma^+\prec\gamma^-$ implies that $\gamma^+$ and $\gamma^-$ are Legendrian isotopic.
\end{Rem}
\begin{Rem}
One might be concerned by the resemblances between our relation and the one given by Arnold in \cite{Arlag2} and \cite{Arlag1}. We briefly recall his definition and emphasize on the differences between this one and ours.

Arnold defined some relations on Lagrangian immersions and Legendrian submanifolds which he called Lagrange and Legendre cobordism. For a manifolds $M$, these relations study exact Lagrangian immersions in $T^*(M\times [0,1])$ (or more generally in $T^*B$ with $B$ with non-empty boundary) which coincide on the reduced slices at $0$ and $1$ with some prescribed exact Lagrangian submanifolds in $T^*M$, and respectively some Legendrian immersion in $Y\times T^*[0,1]$ where $Y$ is a contact manifold. Those relations were symmetric and reflexive, and he defined with them some groups which were later algebraically characterized and computed by Audin in \cite{Audincalculcobordisme}. A variation of those were studied by Ferrand in \cite{Ferrandcob} which also was characterized using a result of Fuchs and Tabachnikov.

Despite their resemblances with the one we propose here, those relations are disjoint from ours. For instance the embedding condition is important in our definition and the symplectic manifold where Lagrangians are studied is different. Even in the case where $Y$ is the $1$--Jet space of a $1$--dimensional manifold $M$ so that there is a symplectomorphism between a compact part of the symplectisation of $J^1(M)$ and $T^*(M\times [0,1])$; this identification necessitates a rescaling of the cotangent fibers by the parameter $t\in [0,1]$. This rescaling is a manifestation of the non-symmetry of the relation and send trivial cylinders in our context to cylinders with different geometrical meaning in the context of \cite{Arlag2} and \cite{Arlag1}.
\end{Rem}
\subsection{Proof of Theorem \ref{Theoremprinc}}\label{isotopieconc}

This subsection is devoted to the proof of \ref{Theoremprinc}  which asserts that the notions of Lagrangian concordance and cobordism are well defined on the set of isotopy classes of Legendrian submanifolds. The main difficulty is the fact that the graph of a Legendrian isotopy is not  a Lagrangian cylinder in general. Conversely, a Lagrangian cylinder which is the graph of a (smooth) isotopy may not be the graph of a Legendrian isotopy.

\begin{proof}[Proof of \ref{Theoremprinc}]
Let $f_t$ be the compactly supported one parameter family of contactomorphisms which realize the isotopy. Denote by $H$ the Hamiltonian on $X\simeq \R\times Y$ whose flow realizes the lifts of $f_t$ as in \ref{symplectisation}.

Now define a new function $H'\co \R\times Y\times [0,1]\appl\R$ such that:
\begin{itemize}
 \item {$H'(v,p,t)=H(p,v,t)$ for $v>T$}
\item {$H'(v,p,t)=0$ for $v<-T$}
\end{itemize}
Denote by $\phi_t$ its Hamiltonian flow.

Then $\phi_1(\R\times K^+)$ coincides with $\R\times K^+$ near $+\infty$ and $\R\times K^-$ near $-\infty$. Since $\R\times K^-$ is Lagrangian and $\phi_1$ is a Hamiltonian diffeomorphism, it follows that $\phi_1(\R\times K^-)$ is a concordance between $K^-$ and $K^+$.

One might worry about the fact that we consider a vector field on a non-compact manifold and hence the flow might not be defined, but in this case outside of a compact set the vector field coincides with one which admits a flow (by its very construction) then so does our new vector field.
\end{proof}
\begin{Rem}
This proof appears to be really simple once \ref{Cont-Ham} is proved. Furthermore we have lots of freedom on the value of the function $H'$, for example we could choose $H'(v,p,t)=\eta_T(t)\cdot H(v,p,t)$ where $\eta_T$ is $0$ on $(-\infty,-T)$ and $1$ on $(T,+\infty)$, increasing on $[-T,T]$ and with small higher order derivatives. For large $T$ this gives, from the point of view of the author, the best concordance for keeping track of the information from the original isotopy. One might be able to compute holomorphic curves with boundary on those cylinders with a conveniently chosen almost complex structure. The naive one, $\phi_1^*(J)$ for $\R$--invariant $J$, is however not admissible for Legendrian contact homology purpose since it is, in general, not admissible at infinity.
\end{Rem}


\section{Properties}\label{behaviourcla}
In this section we discuss the algebraic properties of the Lagrangian concordance. We then discuss the case of Lagrangian surfaces in some fillings with Legendrian boundary and deduce a criterion to compute in some case the slice genus of knots.

\subsection{Behavior of the classical invariants under Lagrangian cobordisms}\label{Behaviour}

\begin{figure}[ht!]
\labellist
\small\hair 2pt
\pinlabel {$f$} [bl] at 173 167
\pinlabel {$\nabla f$} [tl] at 181 11
\endlabellist
\begin{center}
\includegraphics[width=5cm]{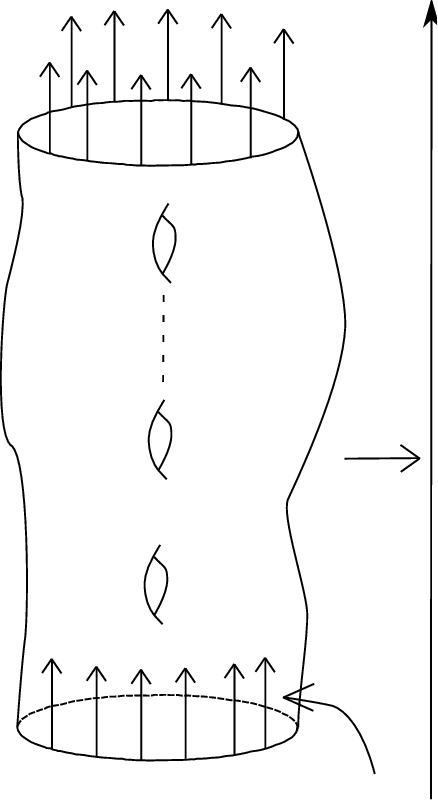}
 \end{center}
 \caption{Height function on $\Sigma$.}
 \label{genreg}
 \end{figure}

The aim of this section is to prove \ref{cobinvariance}.
\begin{proof}[Proof of \ref{cobinvariance}]
 For the invariance of the rotation number just note that $\gamma^-$ and $\gamma^+$ are homologous on $\Sigma$. Since $\Sigma$ is a Lagrangian surface, the loops $l^-(s)$ and $l^+(s)$ of \ref{Clainv} are homologous in $H^1(\Lambda(2))$. Hence we have that $\mu(\gamma^-,S)=\mu(\gamma^+,S\cup[\Sigma])$. The discussion of \ref{Clainv} therefore implies the first part of the theorem.

For the behavior of the Thurston--Bennequin number, recall that multiplication by $J$ will send tangent vectors to a Lagrangian surface to normal ones. On $\Sigma$ take the vector field $\nabla f$ (for the $\difff$--compatible metric) where $f$ is the height function shown in \ref{genreg}.
It has $-\chi(\Sigma)$ zeros of index $-1$ and coincides with $\frac{\del}{\del u^i}$ ($u^i=u^+$ or $u^-$) around the punctures.
Hence any nowhere vanishing vector field $V$ extending $\frac{\del}{\del u^+}$ around the first puncture will have winding number $\chi(\Sigma)$ with respect to $\frac{\del}{\del u^-}$ around the second one.

Since, by definition, around the puncture $\frac{\del}{\del u^i}$ is sent by $dL$ to $\frac{\del}{\del v}$ and since $J\frac{\del}{\del v}$ is $R_\al$ we have that $(JV',W)$ (where $V'$ is the image of $V$ under $dL$ and $W$ is a positive orthogonal to $JV'$ in $N(\Sigma)$) is then a framed cobordism from $(K^-,-tb_{K^-})$ to $(K^+,\sigma)$ where $\sigma$ differs from $-tb_{K^+}$ by $\chi(\Sigma)$ twists. 
\end{proof}

As a corollary, one concludes:
\begin{Cor}\label{invariance}
 If $\gamma^-\prec\gamma^+$ with cylinder $C$ then:
\begin{align*}
 &r(\gamma^-,[S])=r(\gamma^+,[S\cup C])\\
&tb(\gamma^-)=tb(\gamma^+)
\end{align*}
\end{Cor}
 \begin{Rem}
 An examination of Gromov--Lees theorem (Lees \cite{Lees}, Eliashberg and Mishachev \cite[Theorem 16.3.2]{EliMiH}) in this context together with the discussion of \ref{Clainv} gives that the rotation number is the only obstruction to the existence of immersed Lagrangian surfaces in the symplectisation between two Legendrian knots. The behavior of the Thurston--Bennequin invariant is a manifestation of the embedding condition.
\end{Rem}
\begin{Rem}
 \ref{cobinvariance} implies that, unlike in the smooth case, the study of the general cobordism relation is interesting. It will be further investigated by Ekholm, Honda and K\'alm\'an in \cite{KalHoEk}.
\end{Rem}

\subsection{Lagrangian surfaces in Stein filling} $\label{cob}$

The question of Lagrangian concordance naturally leads to the question of whether
or not a Legendrian knot bounds a Lagrangian surface in a symplectic filling.

Let us first study the trivial example which will be relevant in the following:

consider in $D^4\subset\C^2\simeq\Qu$ the disk $D^2=\{z_1=0\}$, which
is a holomorphic disk for the standard K\"ahler structure (namely multiplication
by $i$). If we change the K\"ahler structure by an orthogonal one (eg
multiplication by $j$), this disk turns out to be a Lagrangian one. Now
consider $K_0=\del D ^2\subset S^3$, it is a trivial knot in the $3$--sphere (actually
this is a fiber of the Hopf fibration) which turns out to be Legendrian for the $j$--convex contact structure $\eta_0$. We will argue below that for this knot $$r(K)=0\ \text{and}\ tb(K)=-1=TB(\L_0)$$ where $\L_0$ is the smooth isotopy class of the trivial knot and $TB(\L)$ denote the maximal Thurston--Bennequin number among the Legendrian representatives of a knot type $\L$.

Going back to the general case, through the discussion $X$ will be a Stein
surface with boundary $Y$ endowed with the contact structure $\xi$ induced by
complex tangencies (which is well-known to be tight) and $\L$ will be the
smooth isotopy class of a knot.

Suppose that there is an oriented Lagrangian surface
$L\co \Sigma\emb X$ whose boundary is a Legendrian knot $K\in\L$.
Then Lisca and Mati\'c's adjunction inequality (\cite{LiMa}, and \cite{Rudolph} for $X=D^4$) gives: $$tb(K)+\vert r(K,[\Sigma])\vert\leq
-\chi(\Sigma)=2g(\Sigma)-1$$
Note that $tb(K)$ and $r(K,[\Sigma])$ are numbers associated to $\Sigma$ as in \ref{Clainv}.\\

We can use Gompf's surgery description \cite{Goste} to produce a new Stein
manifold diffeomorphic to $$X(K):=X(K,tb(K)-1)=X\sqcup D^2\times D^2/\{f(S^1\times
D^2)\equiv \mathcal{N}(K)\}$$ where $f$ is the surgery map associated to the framing $tb(K)-1$.

The K\"ahler structure on $X(K)$ is the one induced by $X$ and $D^2\times D^2$ as a subset of $\R^2\oplus
i\R^2$. Now $ D^2\times\{0\}$ is a Lagrangian disk which we use to cap off
$\Sigma$ so that we get a closed Lagrangian surface $\Sigma '$ inside $X(K)$ (notice that it is smooth since the Maslov indices of $K$ on $\Sigma$ and on $D^2$ are equal to zero, see \cite[Lemma 2]{Lasu}). Its self-intersection is consequently: $$\Sigma '\cdot \Sigma '=2g(\Sigma
')-2=2g(\Sigma)-2$$ However the surgery description also tells us that: $$\Sigma '\cdot \Sigma
'=tb(K)-1$$ Combining these formulas and Lisca--Mati\'c's formula we therefore get
$$2g(\Sigma)-1=tb(K)\leq 2g(S)-1 -\vert r(K,[S])\vert$$
where $S$ is any surface in the same homology class as $\Sigma$.

Therefore we get that $g(\Sigma)$ is minimal in its homology class, $tb(K)=2g(\Sigma)-1$ and $r(\Sigma)=0$ which proves \ref{lagsurface}.\hfill\qedsymbol

Intuitively, one can therefore think of Legendrian knots bounding orientable Lagrangian
surfaces as being maximal for the Thurston--Bennequin invariant and, above all,
for the partial ordering given by Lagrangian concordances.
\begin{Rem}
One knows lots of knots yielding examples where Lisca--Mati\'c's inequality is not sharp (eg negative torus knots, connected sum of torus knots...). Consequently, by the previous proposition, such knots cannot
bound orientable Lagrangian surfaces in $D^4$.
\end{Rem}
\begin{Rem}
 Note that knots $K$ such that the trivial knot $K_0$ is Lagrangian cobordant to $K$ are exactly those which bound a Lagrangian surface in $D^4$ since we can patch the Lagrangian disk described at the beginning of the section.
\end{Rem}


\section{Applications and Remarks}$\label{App}$\label{application}

 In this final section we provide some examples of Lagrangian surfaces bounding
Legendrian knots. Those examples lead to some canonical Legendrian representatives of algebraic knots with maximal Thurston--Bennequin number. In the second part of the section, we make a few more remarks about Lagrangian concordances of Legendrian knots, work to be
done in the future by the author.

 \subsection{Algebraic Legendrian knots}\label{algknot}

Let $P\co \C^2\appl\C$ be an irreducible polynomial such that $P(0,0)=0$ and $0$ is a critical
value of $P$ with $(0,0)$ as isolated critical point. The intersection $$K=
P^{-1}(0)\cap S^3_{\epsilon},$$ for small $\eps$, is a
knot (since $P$ is irreducible). Knots arising this way are called {\it algebraic}.

The manifold $S^3_\eps - K$ may be given the structure of a fibration over
$S^1$ via Milnor construction
$$f(z_1,z_2)=\frac{P(z_1,z_2)}{|P(z_1,z_2)|}$$ and the same holds for
$D^4_\eps - P^{-1}(0)$. The fiber of the first fibration is called the Milnor fiber of the singularity and is a
Seifert surface for $K$. One way to formulate the local Thom conjecture is to
say that this Seifert surface is genus minimizing in $D^4$. Notice that the
genus of the fiber of $f$ is given as the genus of $P^{-1}(\delta)\cap D^4_\eps$ for any regular
value $\delta$ of $P$ (see Milnor \cite[Theorem 5.11]{Milcomp}). The aim of this section is to show that we can find a canonical Legendrian
representative of an algebraic knot together with a Lagrangian surface bounded
by it.

On $\C^2\simeq\Qu$ we will consider the K\"ahler structure given by
multiplication by $j$ (instead of $i$) giving the $j$--convex structure on $S^3$. Consider $\delta$ close to $0$ a regular value of $P$ and $P\vert_{S^3_\eps}$ for any $\eps$
sufficiently small. Since the structure induced by $j$ is orthogonal to the one
induced by $i$ and $\Sigma=P^{-1}(\delta)\cap D^4_\eps$ is an $i$--complex curve, $\Sigma$ is
a Lagrangian surface for the symplectic structure we have chosen. However the knot $\Sigma\cap S^3_\eps$ will never be Legendrian for the contact structure induced by the $j$--complex tangencies (except for the trivial polynomial, giving $K_0$). However we have a little more freedom on the contact structure we might choose then. Any Liouville vector field $V$ transverse to $S^3_\eps$ will give a contact form $\al_V=i^*(V\iota\difff)$ whose kernel will be a contact structure isotopic to $\eta_0$ and for which $(D^4_\eps,\difff_0)$ will be a Stein filling. So we are done if we can find such a vector field making $K$ Legendrian.

By Cartan's formula, Liouville vector fields are dual to primitives of $\difff_0 $ and since there is a ``canonical'' primitive to $\difff_0$ (given by $\frac{\del}{\del r}$) we have that Liouville vector fields are dual to exact forms. We denote by $V_f$ the dual to the exact form $df$ (so $\difff(V_f,W)=\theta(W)+df(W)$). We seek then a function $f$ such that:
\begin{itemize}
 \item [1.]{$\difff_0(V_f,\dot{\gamma})=0$}
\item[2.]{$g(V_f,\frac{\del}{\del r}) >0$}
\end{itemize}
Let $\gamma$ be a parametrization of $K$. For $1$ we deduce $\theta(\dot{\gamma})+df(\dot{\gamma})=0$. And so $$f(\gamma(s_0))=-\int_0^{s_0}\theta(\dot{\gamma}(s))ds.$$
This is a well defined function on $K$ since $\theta$ is closed on $\Sigma$ and $K$ is null-homologous on $\Sigma$.\\

We now need to see that we can extend this function to $D^4$ with the only restriction given by $2$. This restriction read as:
$$\theta(R_\al)+df(R_{\al})>0$$
Thus to take such an extension we just need to verify that if $r(t)$ is a trajectory of the Reeb flow starting and finishing on $K$ we have that $\int \theta(\dot{r}(t))dt > f(r(0))-f(r(t_1))=\int_0^{r(t_1)}\theta(\dot{\gamma}(s))ds$. Hence we have to see that:
$$\int_r\theta>\int_{\gamma'}\theta$$
Where $\gamma'$ parametrizes a piece of $K$ between $r(0)$ and $r(t_1)$.

\begin{figure}[ht!]
\vspace{8pt}
\labellist
\small\hair 2pt
\pinlabel {$\gamma'$} [br] at 9 95
\pinlabel {$P^{-1}(0)$} [bl] at 121 95
\pinlabel {$r$} [br] at 0 54
\pinlabel {$L$} [bl] at 178 24
\endlabellist
\begin{center}
\includegraphics[width=10cm]{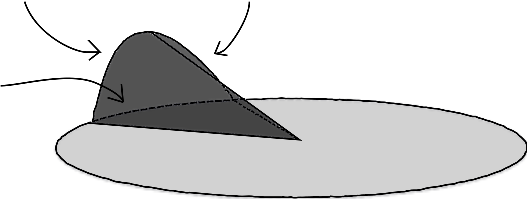}
 \end{center}
 \caption{Cone over a Reeb chord.}
 \label{cone}
 \end{figure}
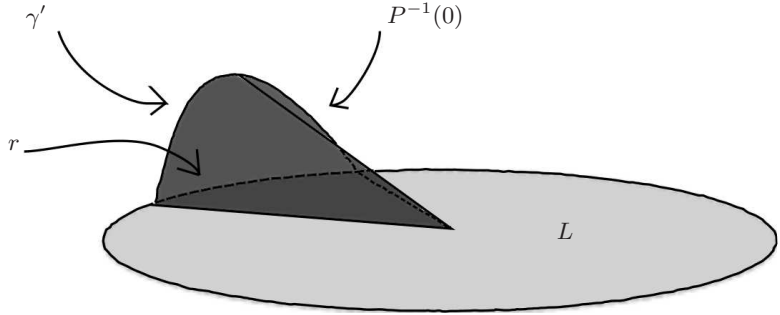

However note that the Reeb flow is given by $j$--complex lines $L$ in $\Qu$ (and so real lines for the complex structure $i$) and that $P^{-1}(0)\cap L$, if not only \{0\}, have to be diffeomorphic to cone over the intersection points $K\cap L$ (because the singularity $P^{-1}(0)$ has the smooth type of the cone over $K$, see \cite[Theorem 2.10]{Milcomp}). So the arc on $K$ starting at $r(0)$ and ending at $r(t_1)$ bounds a piece of disk on $\Sigma$ and the arc $r(t)$ bounds a piece of disk on $L$ and those two disks meet at this cone to form a piecewise smooth disk $D$ (see \ref{cone}). Since one piece of $D$ is symplectic and the other one is Lagrangian we have: $$0<\int_D\difff=\int_r\theta-\int_{\gamma'} \theta$$ Taking $\delta$ small enough does not change this property.
We have then proved that $K_\delta$ is Legendrian for a contact structure isotopic to the standard one which admits a symplectic filling and a Lagrangian surface bounding it. Using \ref{lagsurface}, we have therefore proved:
 \begin{Thm} $\label{Miln}$ Let $K_1$ be an algebraic knot with Milnor fiber
$\Sigma=P^{-1}(\delta)\cap D^4$. Then\\
{\bf (i)}There exists a contact structure on the sphere fillable by $(D^4,\difff_0)$ such that its intersection with $\Sigma$ is a
Legendrian representative $K$ of $\L(K_1)$. \\
{\bf (ii)} The Legendrian representative $K$ satisfies:
\begin{itemize}
 \item{$TB(\L(K))=tb(K)=2g_s(K_1)-1$}
\item{$r(K,\Sigma)=0$}
\end{itemize}
{\bf (iii)} The $4$--ball genus $g_s (K_1)$ of $K_1$ is $g(\Sigma)$.\hfill\qedsymbol\\
\end{Thm}

\begin{Rem}
 Since the contact structure in \ref{Miln} is isotopic to the standard one, one can actually see a Legendrian representative of $\L(K_1)$ directly inside $(S^3,\xi_0)$ with a Lagrangian surface in $D^4$ whose boundary is this representative.
\end{Rem}

 \begin{Rem} The computation of $TB(\L)$ in (ii) is a classical result of
Bennequin, see \cite{Benpf}. Part (iii) is known as the local Thom
conjecture proved by Kronheimer and Mrowka in \cite{KroMr} as we have
shown that the Milnor fiber is genus minimizing among all surfaces in $D^4$
bounded by $K$.
 \end{Rem}

  \subsection{Concluding Remarks}\label{conclusion}

Recall that a topological knot type is {\it Legendrian simple} if the  Legendrian isotopy classes are classified by the Thurston--Bennequin and rotation number. From $\ref{Behaviour}$ and $\ref{Princ}$ we conclude the following corollary.

\begin{Cor} If $\L$ is a Legendrian simple knot type, then any two Legendrian representatives of $\L$ are Lagrangian concordant iff they are Legendrian isotopic.
\end{Cor}

There is also an obvious relation between Lagrangian concordance and Legendrian contact homology. Following  Bourgeois \cite {BouCon} we see that a Lagrangian cylinder between two Legendrian knots could be used to define a map  between the algebras $CH(K^+)$ and $CH(K^-)$ (see Ekholm, Etnyre and Sullivan \cite{EkEtSulegconho}). We, however, will not give a more detailed description of this map for two reasons. First, we have not computed this map for a nontrivial cylinder yet. Moreover Tam\'as K\'alm\'an already gave, in \cite{Kal}, a combinatorial map in contact homology for Legendrian isotopies. So before enlarging our current article, we plan to do two things: give a nontrivial example of such a map, and hopefully be able to compare this map with the one of \cite{Kal} when the cylinder is constructed as in $\ref{Princ}$.

From the result of $\ref{Behaviour}$ we know that on an immersed Lagrangian cylinder between two Legendrian knots, the difference between the Thurston--Bennequin numbers is an obstruction to suppress the double points of the immersion. Suppose that an immersed Lagrangian cylinder is obtained from an embedded cylinder (not Lagrangian). Then the double points of this immersion arise in pairs with opposite signs. To this kind of double points, one can associate another algebraic invariant: the Maslov index of the pair $\mu(x,y)$. We expect this number to be related to the difference of the Thurston--Bennequin numbers. In the case where we actually have an embedded Lagrangian cylinder between two Legendrian knots and perturb it to obtain an immersion with two transverse double points, then the Maslov index of this pair of points is equal to 1 (compare Lalonde \cite{Lasu}). We therefore wish to formulate the following  ``vague'' conjecture:
\begin{Conj} Let $C$ be an immersed Lagrangian cylinder obtained by a deformation of a smooth concordance between two Legendrian knots, and let $\{x_i,y_i\},i\in \{1\cdots k\}$ be cancelable pairs of double points together with $u_i\in\pi_2(x_i,y_i)$ some Whitney disks (these exist by the hypothesis that the pairs are cancelable). Then $$\Sigma_{i=1}^k(\mu(x_i,y_i,u_i)-1)=tb(K^+)-tb(K^-).$$
\end{Conj}

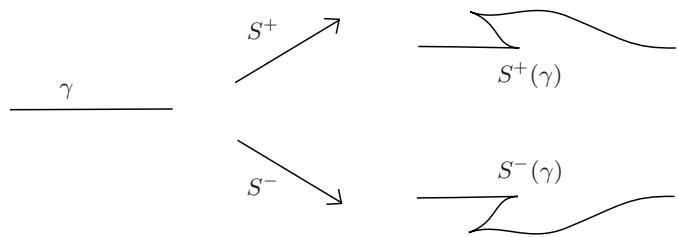
\begin{figure}[ht!]
\labellist
\small\hair 2pt
\pinlabel {$\gamma$} [bl] at 25 71
\pinlabel {$S^+$} [br] at 147 103
\pinlabel {$S^-$} [tr] at 147 34
\pinlabel {$S^+(\gamma)$} [tl] at 260 95
\pinlabel {$S^-(\gamma)$} [bl] at 260 25
\endlabellist
\begin{center}
\includegraphics[height=3cm]{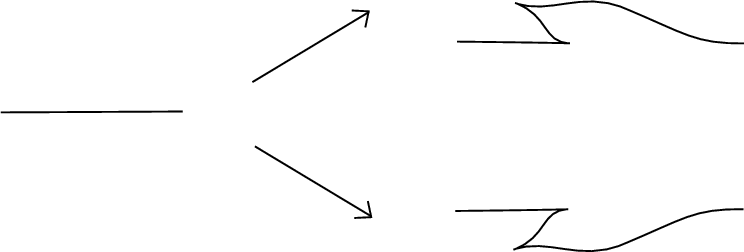}
 \end{center}
 \caption{Front projection of the stabilization of a Legendrian knot.}
 \label{Convexe}
 \end{figure}

We now finish the paper proving one last result concerning the behavior of Lagrangian concordance under stabilization. We recall that the stabilization $S^+$ or $S^-$ of any Legendrian knot $K$ can be defined (see Etnyre and Honda \cite{EtnHoK1}) by the operation which consists of exchanging an arc $\gamma$ of $K$ (in its standard neighborhood) in the way described by \ref{Convexe}. Now suppose that we have a Lagrangian concordance $C$ between $K_0$ and $K_1$ and fix a neighborhood $N$ of $C$ symplectomorphic to the symplectisation of the standard neighborhood of $K_0$. Then replacing $\R\times\gamma$ by the band $\R\times S^{(+,-)}(\gamma)$ gives a Lagrangian concordance $S^{(+,-)}(C)$ between the stabilized knots. Hence:
\begin{Prop} If $C$ be a Lagrangian concordance from $K^-$ to $K^+$ then, $S^{(+,-)}(C)$ is a concordance between $S^{(+,-)}(K^-)$ and $S^{(+,-)}(K^+).$\hfill\qedsymbol
\end{Prop}

\bibliographystyle{plain}
 \bibliography{bibliographie}

\end{document}